\documentclass[12pt,reqno]{amsart}

\usepackage[dvips]{graphicx} 

\usepackage{
hyperref
}

\usepackage{breakurl}   

\usepackage{epigraph}

\theoremstyle{definition}

\numberwithin{equation}{section}

\newcommand\R {{\mathbb R}}

\newcommand\astr{{{}^\ast\hspace{-1pt}\R}}

\renewenvironment{quote}
               {\list{}{}%
                \item\relax}
               {\endlist}

\title{Leibniz's contested infinitesimals: Further depictions}

\author[Mikhail Katz]{Mikhail G. Katz} \address{Department of
  Mathematics, Bar Ilan University, Ramat Gan 5290002 Israel}
\email{katzmik@math.biu.ac.il}

\author[Karl Kuhlemann]{Karl Kuhlemann} \address{Gottfried Wilhelm
  Leibniz University Hannover, D-30167 Hannover, Germany}
\email{kus.kuhlemann@t-online.de}


\begin{document}

\begin{abstract}
We contribute to the lively debate in current scholarship on the
Leibnizian calculus.  In a recent text, Arthur and Rabouin argue that
non-Archimedean continua are incompatible with Leibniz's concepts of
number, quantity and magnitude.

They allege that Leibniz viewed infinitesimals as contradictory, and
claim to deduce such a conclusion from an analysis of the Leibnizian
definition of quantity.  However, their argument is marred by numerous
errors, deliberate omissions, and misrepresentations, stemming in a
number of cases from flawed analyses in their earlier publications.

We defend the thesis, traceable to the classic study by Henk Bos, that
Leibniz used genuine infinitesimals, which he viewed as fictional
mathematical entities (and not merely shorthand for talk about more
ordinary quantities) on par with negatives and imaginaries.
\end{abstract}

\thispagestyle{empty}



\maketitle
\tableofcontents


\epigraph{It is true that, in my view, infinites are not wholes, and
  the infinitely small are not magnitudes.  My metaphysics banishes
  them from its territories, and I don't accord them residence except
  in the imaginary spaces of the geometric calculus, where such
  notions are only as proper as the roots that one calls imaginary.
  The role I played in making this calculus of infinitesimals
  appreciated does not make me sufficiently enamored to push them
  beyond common sense.
\\ --Leibniz to Fontenelle, 9 september 1705 (translation ours).}

\section{Introduction} \label{sec_introduction}

In 2022, Arthur, Rabouin, and others published a brief two-page
broadside \cite{Ar22}, essentially free of serious content, against
the work of scholars who challenge the currently received dogma on the
Leibnizian calculus specifically and other historical
infinitesimalists more generally.  A detailed rebuttal appeared in
Bair et al.~\cite{22a} (a brief response appeared in \cite{23a}).
Apparently feeling that the matter requires a more thorough treatment,
Arthur and Rabouin recently published a 17-page follow-up~\cite{Ar24},
containing a critique of what they describe as ``the vast program
launched in the last ten years or so
{\ldots}\;to reinterpret the history of infinitesimal techniques,
{\ldots}''%
\footnote{Arthur and Rabouin \cite[note 2]{Ar24}.}
We are grateful to Arthur and Rabouin for giving us an opportunity to
develop a detailed assessment of their position.  As we analyze below,
their new article contains significant errors and misinterpretations
of Leibniz's texts.

Richard Arthur and David Rabouin (henceforth AR) announce that they
are pursuing a new strategy so as to tackle the question of the
existence of infinitely small quantities in the Leibnizian calculus.
Such a strategy involves ``the study of his conception of quantity,
number or magnitude.''%
\footnote{Op.~cit., p.\;27.}
They claim to have shown that Leibnizian views on what a
\emph{quantity}, a \emph{magnitude} or a \emph{number} have to be,
``are not compatible with the existence of non-Archimedean instances
thereof.''%
\footnote{Op.~cit., p.\;28.}
They claim further that the ongoing debates about whether nonstandard
analysis (NSA) is an appropriate framework for reinterpreting the
Leibnizian calculus rest on a series of misunderstandings, primarily
on a ``systematic conflation of the question of use with that of
existence of infinitely small quantities.''%
\footnote{Op.~cit., p.\;27.}
(on use and existence see Section~\ref{s8}).

In reality, AR themselves succumb to some significant conflations.  We
point out three major errors in their argumentation:
\begin{enumerate}
\item
the conflation of existence in nature and existence in mathematics,
\item
the conflation of infinite wholes and bounded infinities
(\emph{infinita terminata}),
\item
the conflation of homogeneity and comparability.
\end{enumerate}
In Sections~\ref{s2} through \ref{s4}, we address these errors in more
detail.  Additional misrepresentations and misconceptions in AR's text
are analyzed in Sections~\ref{s5} through \ref{s8}.  Section~\ref{s9}
contrasts the position of AR and that of Henk Bos and other scholars.
Section~\ref{s10} chronicles Rabouin's about-face concerning the
presence of ``limit-processes'' in the Leibnizian calculus.

\section{Existence in nature and existence in mathematics}
\label{s2}

AR's main argument for the incompatibility of Leibniz's views with
non-Archimedean continua is that Leibniz described his inassignables
as \emph{impossible entities}, as \emph{contradictory} and as
\emph{absurd}.  According to AR, Leibniz proved the absurdity of
infinitely small lines.  They claim further: ``we can exhibit a text
in which Leibniz finally elaborated this proof in full detail: the
notes he took on Froidmont.''%
\footnote{Op.~cit., p.\;38.}
AR go on to present the following Leibnizian passage:

\begin{quote}
``Therefore at some time the ratio of what is traversed to the
  unassignable whole must necessarily begin to be assignable, but such
  a point or instant in which that transition comes about is
  impossible. Therefore such a \emph{motion}, and so such a line, must
  also be held to be impossible, and there is no infinite yet bounded
  straight line. Whence it also follows that infinitely small straight
  lines are chimeras, although nevertheless useful, like imaginary
  roots \dots''%
\footnote{Leibniz as translated in op.~cit., p.\;39; emphasis added.}
\end{quote}
Note that Leibniz speaks of motion (a \emph{physical} notion).  In
fact, a few lines earlier AR acknowledge that
\begin{quote}
The general context is a part of the notes in which Leibniz goes back
to the possibility of interpreting physical \emph{conatus} as
infinitely small lines traversed in a moment.%
\footnote{Op.~cit., pp.\;38--39.}
\end{quote}
Evidently, the issue is the question of the existence of
infinitesimals in the physical realm.  Failing to distinguish the
latter sufficiently from the mathematical realm,%
\footnote{Arthur has a history of conflating these two Leibnizian
realms by applying actuality to the mathematical realm, as we discuss
in Section~\ref{s21b}.}
AR argue as follows:
\begin{quote}
There is no ambiguity about the conclusion here: Leibniz considered
such a reasoning as establishing the \emph{impossibility} of infinite
bounded lines, so that they should be ranked with infinitely small
lines in the store of \emph{chimerical} mathematical entities.  But
the most important fact is to note that from 1676 to the reflections
on Froidmont, we see Leibniz considering infinitesimal lines as
\emph{absurd}, even when not taken as genuine magnitudes.%
\footnote{Arthur and Rabouin \cite[p.\;39]{Ar24}; emphasis added.}
\end{quote}
So far, infinitesimals are \emph{impossible}, \emph{chimerical}, and
\emph{absurd}.  But then AR jump to the following conclusion:
\begin{quote}
As we have established, he attempted regularly to prove that they
involve a \emph{contradiction}.%
\footnote{Ibid.; emphasis added.}
\end{quote}
Did Leibniz attempt to prove that infinitesimals involve a
contradiction, as AR claim?  As AR correctly point out some pages
earlier,
\begin{quote}
According to Leibniz, mathematical existence can be equated with
non-contradiction.%
\footnote{Op.~cit., p.\;27, note 6.}

\end{quote}
If so, how could Leibniz have used infinitesimals if they are
contradictory according to AR?  We will examine the matter in
Section~\ref{s21}.

\subsection{Accidental and absolute impossibility}
\label{s21}

What AR appear to ignore is the fact that Leibniz distinguished
between existence in the phenomenal realm (i.e., physical realm) and
mathematical existence (which also includes fictional entities).

For Leibniz, \emph{fictional} did not necessarily mean
\emph{contradictory}.  Crucial here is Leibniz's distinction between
absolute impossibility and accidental impossibility.  Already in the
\emph{Confessio Philosophi},%
\footnote{Leibniz \cite{Le72}, 1672/3, p.\,128.}
Leibniz speaks of the concept of ``impossible by accident'' and
contrasts it with ``abso\-lute impossibility'' i.e., contradiction.
He gives the examples of a species with an odd number of feet, and an
immortal mindless creature, which are, according to him,
\emph{harmoniae rerum adversa} i.e., ``contrary to the harmony of
things,''%
\footnote{Leibniz as translated by Sleigh in \cite[p.\;57]{Le04c}.}
but \emph{not} necessarily contradictory.

In his 1683 text \emph{Elementa nova matheseos universalis}, Leibniz
explains that some mathematical operations cannot be performed in
actuality, but nonetheless one can exhibit
\begin{quote}
``a construction in our characters'' (\emph{in nostris
  characteribus}),%
\footnote{Leibniz \cite[p.\;520]{Le83}.}
\end{quote}
meaning that one can carry out a formal calculation, such as those
with imaginary roots, regardless of whether the mathematical notions
involved idealize anything in nature.  Leibniz goes on to discuss in
detail the cases of imaginary roots and infinitesimals, which to him
constitute examples of accidental impossibility.%
\footnote{For a detailed textual analysis, see
\cite[Section~5.3]{24a}.  The distinction was also analyzed by Eklund
\cite{Ek20}.}

Something contradictory is absolutely impossible and can neither occur
in nature nor be allowed in mathematics.  It may be absurd to assume
the existence of negatives, imaginaries, or infinitesimal quantities
in nature.%
\footnote{AR's article contains several relevant passages from
Leibniz: ``whether the nature of things tolerates quantities of this
kind is for the metaphysician to discuss; it is enough for the
geometrician to demonstrate what follows from supposing them''
\cite[p.\;35]{Ar24}; ``I would indeed admit these infinitely small
spaces and times in geometry for the sake of discovery, even if they
are imaginary.  But I doubt whether they can be admitted in nature''
(op.~cit., p.\;37).}
However, such entities are not contradictory, but useful for
mathematics.  Leibniz therefore described them as well-founded
fictions.  Such a viewpoint involves no conflation of use and
existence of infinitesimals; see further in Section~\ref{s8}.  AR's
conclusion that
\begin{quote}
from 1676 to the reflections on Froidmont, we see Leibniz considering
infinitesimal lines as \emph{absurd},
{\ldots}~
he attempted regularly to prove that they involve a
\emph{contradiction}%
\footnote{Arthur and Rabouin \cite[p.\;39]{Ar24}; emphasis added.}
\end{quote}
is a non-sequitur, because infinitesimals can be ``absurd'' without
being contradictory.  AR claim further that
\begin{quote}
Leibniz's metaphysics claims that only simple substances should be
called `real' \emph{stricto sensu} and it is not easy to understand in
this framework the distinction between fictitious and real entities in
mathematics.%
\footnote{Op.~cit., note~46.}
\end{quote}
In fact, AR exaggerate the difficulty involved.  Certainly, Leibniz
considered only entities in a metaphysical realm as truly \emph{real}.
But as far as mathematics is concerned, ordinary mathematical entities
are those that were familiar to 17th century scholars from their study
of the Greeks.  Meanwhile, entities of a new type: imaginaries,
negatives, and infinitesimals%
\footnote{Leibniz specifically compares the status of negatives and
infinitesimals: ``[J]ust as I deny the reality of a ratio one of whose
terms is a quantity less than nothing, I also deny that properly
speaking there exists an infinite or infinitely small number, or an
infinite or infinitely small line, {\ldots}'' This passage is
translated in \cite[p.\;30]{Ar24} by AR, who fail to examine its
implications.}
would be \emph{fictions}; see~\cite{24e} for a fuller discussion.

\subsection{An epigraph and an artful ellipsis}

AR's epigraph furnishes an interesting case study.  They quote Leibniz
as follows:
\begin{quote}
``Il est vrai que chez moi les infinis ne sont pas des touts et les
  infiniment petits ne sont pas des grandeurs. (\ldots)%
\footnote{The ellipsis is in the original.}
La part que j’ai eue \`a faire valoir ce calcul des infinit\'esimales
ne me rend pas assez amoureux pour les pousser au-del\`a du bon\;sens.''
(Leibniz \`a Fontenelle, 9 septembre 1705)%
\footnote{Leibniz as quoted by Arthur and Rabouin
\cite[Epigraph]{Ar24}.  A translation of the full passage appears in
note~\ref{f22}.}
\end{quote}
We note the ellipsis in the middle of the quoted passage.  The passage
as presented makes it appear as though Leibniz is telling Fontenelle
that genuine infinitesimals would go beyond common sense (``bon
sens'').  Was that truly Leibniz's intention?  A fuller quotation is
readily found in an article by Bella.  Here is the full passage, as
quoted by Bella:
\begin{quote}
``Il est vrai que chez moi les infinis ne sont pas des touts et les
  infiniment petits ne sont pas des grandeurs.  \emph{Ma
  m\'etaphysique les bannit de ses terres, et je ne leur donne retrait
  que dans les espaces imaginaires du calcul g\'eom\'etrique, ou ces
  notions ne sont de mise que comme les racines qu’on appelle
  imaginaires}.  La part que j’ai eue \`a faire valoir ce calcul des
  infinit\'esimales ne me rend pas assez amoureux pour les pousser
  au-delà du bon\;sens.''%
\footnote{\label{f22}Leibniz as cited in Bella \cite[notes 77,
  86]{Be14}; emphasis added.  Translation: ``It is true that, in my
view, infinites are not wholes, and the infinitely small are not
magnitudes.  My metaphysics banishes them from its territories, and I
don't accord them residence except in the imaginary spaces of the
geometric calculus, where such notions are only as proper as the roots
that one calls imaginary.  The role I played in making this calculus
of infinitesimals appreciated does not make me sufficiently enamored
to push them beyond common sense.''}
\end{quote}
It emerges that, contrary to the picture AR seek to paint, the meaning
of the Leibnizian passage is as follows:

\begin{enumerate}
\item
it is Leibniz's metaphysics that expels infinitesimals from ``ses
terres'';
\item
Leibniz nonetheless leaves room for infinitesimals in the imaginary
spaces of geometric calculation;
\item
he compares infinitesimals to imaginary roots.
\end{enumerate}
Accordingly, genuine infinitesimals are fine in mathematics just as
imaginary roots are, but insisting on their existence beyond the realm
of mathematics would push them beyond common sense, according to
Leibniz.  By suppressing the crucial middle sentence, AR changed the
apparent meaning of the Leibnizian passage.  The letter to Fontenelle
proves precisely the opposite of what AR make it out to say by means
of a skillful use of an ellipsis.%
\footnote{Earlier skillful uses of the truncation technique to score
an apparently syncategorematic point in the writings of Ferraro
\cite[p.\;29]{Fe08}, Goldenbaum \cite[p.\;76]{Go08a}, Ishiguro
\cite{Is90}, and Tho \cite[p.\;71]{Th12} were analyzed in Bascelli et
al.~\cite{16a}.}

\section{Infinite wholes and bounded infinities}
\label{s3}

Most commentators agree that Leibniz rejected infinite wholes as
contradictory.  From this, AR seek to deduce that he therefore also
rejected the notion of a bounded infinite line (\emph{infinitum
terminatum}),%
\footnote{The notion of \emph{infinitum terminatum} is analyzed in
Section~\ref{s7}.}
because the latter would consist of an infinite number
of finite segments and such a ``number'' would have to be an infinite
whole.

There is indeed a subtle difference between an infinite
\emph{multitude} of things (e.g., segments) as a whole, and a bounded
but infinitely extended line consisting of an infinite \emph{number}
of segments.  It is plain that infinite multitude was not a number for
Leibniz, but the infinite extension of these segments was a tool
useful in geometry.  Leibniz distinguished between bounded infinity
and infinite wholes, and indeed contrasted them.%
\footnote{\label{f23}Thus, referring to \emph{infinita terminata} as
simply `infinita', Leibniz writes: ``quantitas interminata differt ab
infinita'' \cite{Le23} A.VII. 6. 549.}
As noted by Knobloch, the distinction between \emph{infinita
terminata} and \emph{infinita interminata} was elaborated in
Proposition\;11 of the treatise \emph{De Quadratura Arithmetica}
(DQA):
\begin{quote}
[Leibniz] distinguished between two infinites, the bounded infinite
straight line, the recta infinita terminata, and the unbounded
infinite straight line, the recta infinita interminata.  He
investigated this distinction in several studies from the year 1676.
Only the first kind of straight lines can be used in mathematics, as
he underlined in his proof of theorem 11 [i.e.,
  \emph{Propositio~XI}].%
%
%
\footnote{Knobloch \cite[p.\;97]{Kn99}.}
\end{quote}
Leibniz mentioned the distinction in a 29\;july\,1698 letter to
Bernoulli.  Here Leibniz analyzes a geometric problem involving
unbounded infinite areas and states an apparent paradox: ``Therefore
the two infinite spaces are equal, and the part is equal to the whole,
which is impossible.''%
\footnote{``Ergo haec duo spatia infinita sibi sunt aequalia, pars
toti, quod est absurdum'' (Leibniz \cite[p.\;523]{Le98}).}
To resolve the paradox (i.e., the clash with the part-whole
principle), Leibniz exploits the notion of bounded infinity.  Far from
\emph{creating} paradox, bounded infinities serve to \emph{resolve}
one.%
\footnote{For a detailed textual analysis, see \cite[Section~1]{23h}.}
Further details on \emph{infinita terminata} appear in
Section~\ref{s7}.

\subsection{Whose \emph{wholes}?}
\label{s31}

In their Conclusion, AR claim to ``have shown that a quantity or
magnitude is always defined as a `whole' endowed with a certain number
of (homogeneous) parts''%
\footnote{Arthur and Rabouin \cite[p.\;41]{Ar24}.}
and go on to allege that since Leibniz rejected infinite wholes, he
would necessarily also reject infinite quantity.  How did they show
this exactly?  Their section~3 entitled ``Number and quantity
according to Leibniz" quotes Leibniz as follows:
\begin{quote}
``Quantity is that which belongs to a thing insofar as it has all its
  parts, or, on account of which it is said to be equal to, greater
  than, or less than, another thing (with each homogeneous with the
  other), or can be compared with it.''%
\footnote{Leibniz as translated in Arthur and Rabouin
\cite[p.\;29]{Ar24}.}
\end{quote}
Note that `wholes' are not mentioned in Leibniz's passage.
Immediately following this quote, AR proceed to paraphrase Leibniz:
\begin{quote}
According to this definition, a quantity is a \emph{whole} constituted
of parts and only things to which the part-whole axiom can be applied
can count as quantities \emph{stricto sensu}.%
\footnote{Ibid.; emphasis on `whole' added.}
\end{quote}
Note that the term \emph{whole} only appears in AR's paraphrase.
Leibniz does not use the term in his definition of \emph{quantity}.  A
similar paraphrase occurs when AR present a Leibnizian definition of
magnitude:
\begin{quote}
Leibniz always defines continuous magnitudes, exactly like quantities,
in terms of whole and parts: \emph{``Magnitude is that which is
designated by a number of congruent parts''}.  If infinite
\emph{wholes} are contradictory, infinite magnitudes are contradictory
notions from their very definition, and so are infinitely small parts
of any finite magnitude.%
\footnote{Arthur and Rabouin \cite[p.\;30]{Ar24}; emphasis added.}
\end{quote}
In this passage as well, the term \emph{whole} is mentioned not by
Leibniz but by AR, whose ``contradictory'' conclusion therefore
depends on their paraphrase.

It emerges that AR have shown only that their \emph{paraphrases} of
the Leibnizian definitions of \emph{quantity} and \emph{magnitude} may
lead to a rejection of infinite instances thereof, rather than the
definitions as they appear in Leibniz.

\subsection{Do infinite numbers require infinite sets?}
\label{s32}

AR claim that
\begin{quote}
Leibniz's views about quantities, numbers and magnitude are not
consistent with the NSA interpretation.  This inconsistency takes the
simplest form: the very notion of ``quantity'' in Leibniz, be they
numbers or magnitudes, is \emph{not compatible} with the fact of being
``non-Archimedean''.%
\footnote{Op.~cit., p.\;28; emphasis added.}
\end{quote}
Such an alleged `incompatibility' is due to AR's belief that an
infinite quantity would necessarily be tied up with infinite sets as
per AR (see Section~\ref{s31}), and therefore infinite quantities must
be as contradictory as infinite sets in Leibniz's world.  AR detail
their objections as follows:
\begin{quote}
There are several versions of nonstandard analyses [sic] depending on
the kind of logical tools which are employed in order to introduce
non-Archimedean quantities in mathematical analysis.  Most of them,
however, suppose strong set theoretic commitments like the existence
of ultrapowers or conservative extensions of ZFC (as in various
versions of Internal Set Theory).  These commitments are not
compatible with Leibniz's conviction that any \emph{infinite totality}
is a contradictory notion {\ldots}%
\footnote{Op.~cit., note 5; emphasis added.}
\end{quote}
Here AR speak of `strong set theoretic commitments', `existence of
ultrapowers', and `ZFC'.  The passage makes it appear as if infinite
numbers in modern mathematics necessarily require `infinite
totalities' as well as the axiom of choice (part of ZFC), and are thus
inapplicable to interpreting Leibniz, who rejected infinite
totalities.  But AR's modern set-theoric excursion misses the point,
as we explain in Section~\ref{s33}.

\subsection{Conservativity results by Friedman and others}
\label{s33}

In the late 1960s, Harvey Friedman proved the conservativity of not
only~${}^\ast$PA over PA, but the more powerful theory, namely
${}^\ast$PA plus the induction scheme referred to as
$^*\Pi_{\infty}\mathsf{-Induction}$.  Friedman's work was mentioned by
Kreisel.%
\footnote{See Kreisel \cite[p.\;94]{Fr69}.  See also Feferman
\cite[p.\;967]{Fe77} (second theorem in Section~8.7) and Kossak
\cite{Ko22}.}
Friedman's result was followed up by Keisler et al.~\cite{He84}.
Thus, Robinson's question concerning axiomatisations was answered by
Friedman and led to a mature body of work.%
\footnote{See Enayat \cite{En24} for further details.}
%


Here PA denotes Peano Arithmetic, an axiomatic system that does
\emph{not} prove the existence of an infinite set.  The extension
${}^\ast$PA incorporates the standardness predicate (as well as other
axioms elaborated there) and therefore postulates the existence of
nonstandard integers.  By conservativity,~${}^\ast$PA does not prove
the existence of an infinite set, either.  Therefore the existence of
an infinite -- or more precisely, \emph{unlimited} -- integer does not
require the existence of infinite sets.

To emphasize the point, one could refer to infinite integers in
nonstandard models of arithmetic as \emph{ringinals}, because they are
elements of a (semi)ring, to contrast them with infinite numbers of
Cantorian set theory such as ordinals or cardinals (which do not
satisfy the properties of a ring or semiring).  An infinite ringinal
can be thought of as a formalisation of the Leibnizian \emph{infinitum
terminatum}.


AR attempt to argue from ``first philosophical principles'' that an
infinite quantity would necessarily rely on the existence of an
infinite whole, thereby contradicting a basic premise in Leibniz.
However, the result just mentioned undermines their argument in a
mathematically demonstrable way: infinite ringinals are compatible
with Euclid's part-whole axiom, and do not require infinite numbers
\`a la Cantor.

\subsection{``Bias against NSA''?}
\label{s34}

Historians who received their
training in mathematics in the traditional Weierstrassian paradigm
based on naive Cantorian set theory appear to have difficulty
understanding the distinction between unlimited integers, on the one
hand, and Cantorian infinities (ordinals and cardinals), on the other.
To give an example, AR misunderstood the following passage in Katz et
al.~\cite{22b}:
\begin{quote}
[One] notes that Leibniz stresses the distinction between infinite
cardinality and infinite quantity (reciprocal of infinitesimals) [and]
argues that the exchange with Bernoulli precisely refutes [the]
attempt to blend infinite cardinality and infinite quantity so as to
deduce the inconsistency of infinitesimals.%
\footnote{Katz et al.~\cite[p.\;264]{22b}.}
\end{quote}
AR have misread the above passage entirely.  They write:
\begin{quote}
Katz, Kuhlemann, Sherry, Ugaglia and van Atten (2021) maintain, for
example, that Leibniz’s argument against Bernoulli is only an argument
against infinitesimals interpreted as reciprocals of infinite
multiplicities; if there is no number corresponding to an infinite
multiplicity, then there can be no reciprocal of such a number
either.%
\footnote{Arthur and Rabouin \cite[note 27]{Ar24}.}
\end{quote}
So far, AR's summary is accurate.  They continue as follows, still
referring to \cite{22b}:
\begin{quote}
This, \emph{they claim} [i.e., Katz et al.~claim], leaves open the
question of \emph{continuous} quantities (magnitudes), and
infinitesimals as their reciprocals. Notice that this argument is
clearly \emph{an invention} made in order to save their interpretation
from criticism.  Until 2016, Katz and co-authors regularly presented
infinitesimal magnitudes as a kind of numbers (see the conclusion of
Bascelli et al. 2016, 861) and never mentioned the distinction between
\emph{numbers and magnitudes} as a crucial one (see, e.g., Bair et
al., 2013, 896).%
\footnote{Ibid.; emphasis on `they claim', `continuous', `an
invention', and `numbers and magnitudes' added.}
\end{quote}
However, the authors of \cite{22b} claimed no such thing.  The
distinction is not between discrete and continuous quantity or
magnitude, contrary to what is claimed by AR.  The distinction is not
between number and magnitude either, contrary to what they claim.
Rather, the distinction is the Leibnizian distinction between
magnitude and multitude;%
\footnote{\label{f40}For details see e.g., \cite[Section~4.1]{24a}.
}
in other words, between magnitude or quantity (as formalized by
unlimited numbers in modern theories e.g.,~${}^\ast$PA, as discussed
in Section~\ref{s33}), on the one hand, and infinite wholes formalized
as Cantorian infinities (whether ordinal or cardinal), on the other.%
\footnote{\label{f32}Leibniz has a more general notion of an `infinite
whole' that includes unbounded lines and unbounded regions.  See the
beginning of Section~\ref{s3} and sources given there for a discussion
of the paradox of unbounded regions in an exchange of letters with
Bernoulli.}
If there is any \emph{invention} involved, it is due to AR rather than
the authors of~\cite{22b}.  Such a misunderstanding and an attendant
conflation are a persistent feature of AR's writing on Leibniz, as
when they write:
\begin{quote}
[Leibniz] stumbles upon the \emph{impossibility} of infinite wholes as
early as 1672, although he explicitly rejects actual infinitesimals
only in 1676.%
\footnote{Arthur and Rabouin \cite[p.\;36]{Ar24}; emphasis added.}
\end{quote}
But infinite wholes are \emph{multitudes}, whereas infinitesimals are
\emph{magnitudes}.  The former are not merely \emph{impossible} but
actually contradictory according to Leibniz; the latter are not
contradictory.

In this connection, we would like to comment on AR's claim that
\begin{quote}
[T]he matter has often taken a controversial twist after Robinson, any
objection to the faithfulness of this reading being seen as relying on
an ideological \emph{bias\;against\;NSA} in its various forms.%
\footnote{Op.~cit., p.\;27; emphasis added.}
\end{quote}
Recent scholarship (e.g., \cite{21a}) has argued that the Leibnizian
concept of \emph{infinitum terminatum} constitutes a non-Archimedean
phenomenon, which can be formalized by unlimited numbers in modern NSA
or in weaker theories of Peano Arithmetic.  The issue is not
allegations of ``bias against NSA'' (as AR claim) but rather the fact
that non-Archimedean conceptions are not sufficiently taken into
account by AR in their analysis of Leibniz's mature views.  AR
announce a follow-up in~\cite{Ar24b}; it remains to be seen whether it
is more coherent than their~\cite{Ar24}.


\section{Homogeneity and comparability}
\label{s4}

AR present Leibniz's definition of \emph{number} from GM VII, 31:
\begin{quote}
``\emph{Number} is homogeneous with unity, and so it can be compared
  with unity by adding to it or subtracting from it. And it is either
  an aggregate of unities which is called an \emph{integer},
  {\dots}\;or an aggregate of aliquot parts of unity, which is called
  a \emph{fraction}\ldots''%
%
\footnote{Leibniz as translated in op.~cit., p.\;30.}
\end{quote}
They conclude:
\begin{quote}
But this definition of number is incompatible with treating an
infinitesimal as a number in relation to 1 and ``infinitesimals'' are
never included by Leibniz in his lists of numbers
%
{\ldots}\;\emph{Infinitesimals cannot be numbers, according to
Leibniz, since they are not homogeneous with 1.} They cannot be
``compared with unity.''%
\footnote{Op.~cit., p.\;31; emphasis in the original.}
\end{quote}
Here AR took the Leibnizian remarks out of their context.  Leibniz's
definition of \emph{number} does contain the term `compared', but in
the sense of being able add and subtract them (i.e., of being
homogeneous), not in the sense of Euclid's Definition V.4.%
\footnote{According to Euclid Book V, Definition 4, two homogeneous
quantities are comparable if one, multiplied by a finite number, can
exceed the other.}
Ordinary and infinitesimal numbers can indeed be added and subtracted
together, e.g.,~$2+dx$ or~$2-dx$.

Note that this Leibnizian definition mentions both \emph{number} and
\emph{fraction}.  By AR's logic, \emph{fractions}
would similarly be incompatible with infinitesimals.  If so, how would
AR account for Leibniz's definition of infinitesimal
as a \emph{fraction}, namely ``[an] infinitely small fraction, or one
whose denominator is an infinite number''%
\footnote{``fraction infiniment petite, ou dont le denominateur soit
un nombre infini'' Leibniz \cite[p.\;93]{Le02}.}
?

\subsection{Leibniz on Euclid V.4 in 1695}

Regrettably, AR keep ignoring the fact that in 1695 Leibniz made it
clear in two separate texts that infinitesimals violate Euclid's
Definition V.4.%
\footnote{\label{f5}Leibniz \cite[p.\;288]{Le95a},
\cite[p.\;322]{Le95b}.  For a detailed analysis see \cite[Section 1.5,
  pp.\;177--179]{21a}.}
Neither of these Leibnizian texts is even cited in \cite{Ar24}, even
though AR appear to acknowledge that the presence of inassignables
would signal non-Archimedean behavior; see Section~\ref{s54b}.
Instead, they offer
non-sequiturs, such as the following:
\begin{quote}
In the framework of the geometry of the indivisibles, as in the
framework of differential calculus, Leibniz’s view is
\emph{very~clear}: an infinitesimal~$dx$ stands not for a fixed
infinitely small quantity but for a finite variable quantity, which
one can take as small as one wishes (and the same will be true for
infinitely large quantities considered as variable finite quantities
taken as great as one wishes).%
\footnote{Arthur and Rabouin \cite[p.\,35]{Ar24}; emphasis added.}
\end{quote}
In reality, such an interpretation is far from being `very clear'.
AR's claim is merely a restatement of Ishiguro's position, which
remains as unfounded as it was in 1990, as we discuss in
Section~\ref{f29}.  In their Conclusion, AR claim further:
\begin{quote}
Ever since the publication of the Leibnizian differential algorithm,
controversies have occurred amongst mathematicians to determine
whether or not it commits its users to accept what are now called
non-Archimedean quantities.  Leibniz’s own view on that matter, often
repeated, was that his algorithm was \emph{neutral with regard to these
questions}.%
\footnote{Op.~cit., p.\;40; emphasis added.}
\end{quote}
But Leibniz's view was not ``neutral with regard to these questions''
as he made it clear in a pair of 1695 texts that infinitesimals
violate Euclid's Definition V.4.%
\footnote{See note \ref{f5} for the references.}
Leibniz believed that Archimedean paraphrases are possible, but they
amount merely to a reassuring alternative; see further in
Section~\ref{s8}.

\subsection{Ishiguro on referring and designating}
\label{f29}

The syncategorematic interpretation of Leibnizian infinitesimals goes
back to Ishiguro \cite[Chapter~5]{Is90}.  It is the idea that talk
about infinitesimals is just stenography for a more laborious argument
in a purely Archimedean context, \`a la exhaustion arguments.  The
proponents of Ishiguro's interpretation fiercely reject the idea that
mature Leibniz had genuine infinitesimals in mind (i.e., magnitudes in
a non-Archimedean context).

Both sides (let's name them Alice and Bob following \cite{22b}) agree
that on Leibniz's mature view, infinitesimals were \emph{fictions}.
However, Alice argues that \emph{fiction} means stenography as above,
whereas Bob argues that \emph{fiction} means that, unlike his allies
l'Hospital and Bernoulli, Leibniz was not convinced that such
magnitudes could be found in the physical realm.  Being a philosopher,
Ishiguro emphasizes the idea that the term \emph{infinitesimal} does
not ``refer'':
\begin{quote}
Fictions {\ldots}~are not entities to which we refer.  They are not
abstract entities.  They are \emph{correlates of ways of speaking}
which can be reduced to talk about more standard kinds of
entities. Reference to mathematical fictions can be paraphrased into
talk about standard mathematical entities.%
\footnote{Ishiguro \cite[p.\,100]{Is90}; emphasis added.  Moreover,
``The word `infinitesimal' does not designate a special kind of
magnitude.  It does not designate at all'' (Op.~cit., p.\;83).}
\end{quote}
Since such a use of the verb \emph{refer} may not be familiar to our
reader, we note that ``correlates of ways of speaking'' means that
talk about infinitesimals is shorthand/stenography for more laborious
arguments, etc.  A detailed rebuttal of Ishiguro's reading appeared
in~\cite{16a}.  For our purposes, it will suffice to examine
Ishiguro's purported source for the syncategorematic interpretation,
mentioned early in her Chapter 5:
\begin{quote}
Leibniz {\ldots}~maintained that one can have a rigorous language of
infinity and infinitesimal while at the same time considering these
expressions as being \emph{syncategorematic} (in the sense of the
Scholastics), i.e., regarding the words as not designating entities
but as being well defined in the proposition in which they occur.
(Letter to Varignon, 2 February 1702, G. Math. 4 93).%
\footnote{Ishiguro \cite[p.\;82]{Is90}; emphasis added.}
\end{quote}
The parenthetical reference is to Leibniz \cite[p.\;93]{Le02}.  In her
Chapter~5, Ishiguro seeks to create the impression that the letter to
Varignon provides evidence for her interpretation.  Does Leibniz apply
the term \emph{syncategorematic} to infinite magnitudes in this
letter?  We will examine the issue in Section~\ref{s43}.

\subsection{Conflation of physical and mathematical realms}
\label{s21b}

Arthur has a history of conflating Leibniz's physical and mathematical
realms, and applying actuality to the mathematical realm.  Thus, he
claims the following against three Leibniz scholars:
\begin{quote}
(contra Breger, Bosinelli and Antognazza)%
\footnote{The reference is to Antognazza \cite{An15}.}
there is \emph{not} one kind of infinity that applies to mathematical
entities, and another that applies to actual things.%
\footnote{Arthur \cite[p.\,158]{Ar18}; emphasis added.}
\end{quote}
Furthermore, ``Leibniz's actual infinite, understood
syncategorematically,%
\footnote{For an analysis of the syncategorematic interpretation see
Section~\ref{f29}.}
applies to any entities that are actually infinite in multitude,
whether numbers, actual parts of matter, or monads.''%
\footnote{Op.~cit., Abstract.  There are additional passages of this
sort: ``Leibniz's actual infinite, understood syncategorematically,
applies to numbers and to constituents of the actual world in exactly
the same way'' (Op.~cit., p.\,157); ``the actual infinite pertaining to
the created world, insofar as it is an infinite in multitude, is
syncategorematic in exactly the same sense as the actual infinite that
applies to numbers'' (Op.~cit., p.\,158).}

After quoting a Leibnizian passage on the actually infinite division
of \emph{bodies}, Arthur comments: ``This is exactly the same notion
of the actual infinite in multitude as occurs in his mathematics.''%
\footnote{Op.~cit., p.\,166.}
But \emph{bodies} are not the subject matter of mathematics.

On the next page, we are relieved to learn that ``Antognazza
{\ldots}\;insists that the \emph{actual} infinite cannot apply to
\emph{ideal} entities like numbers''%
\footnote{Op.~cit., p.\,167, note 16.}
and we concur with Antognazza.

By page 167 (twelve pages into his article), we locate the source of
Arthur's \emph{syncategorematic actual} interpretation for
\emph{numbers}, namely Ishiguro:
\begin{quote}
Thus as Hid\'e Ishiguro has affirmed, for Leibniz ``there are
infinitely many substances and infinitely many \emph{numbers} in his
sense of infinitely many, i.e., actually more than any finite number
of them, and not merely potentially more'' (Ishiguro 1990, 80).%
\footnote{Op.~cit., p.\,167; emphasis added.}
\end{quote}
Attached to this sentence on page~80 in Ishiguro's text is a
footnote~2.  The footnote provides her two sources.  The sources are
passages from (a) Leibniz's letter to Foucher from 1692/93, and (b)
Leibniz's \emph{New Essays}.  The passages deal with division of
`matter', `creatures', and `things'.  They say nothing about `numbers'
and therefore don't support her (and subsequently Arthur's)
application of \emph{actuality} to numbers.  Does Arthur provide any
primary sources for his application of `syncategorematic actuality' to
numbers?  We examine the issue in Section~\ref{s44}.

Arthur writes furthermore that ``the syncategorematic infinitesimal is
not merely an indeterminate finite part: it is a fictional
\emph{actually} infinitely small element, whose use may be justified
by exploiting the Archimedean property of the continuum, that is, by
substituting for it sufficiently many sufficiently small parts to
demonstrate that the error is null.''%
\footnote{Op.~cit., p.\,170; emphasis added.}

Finally, ``The actually infinite division of actual matter,
syncategorematically understood, issues in \emph{infinitely small
parts} that are not so small that there are not smaller (finite)
parts.''%
\footnote{Op.~cit., p.\,177; emphasis added.}

In conclusion, Arthur's talk about syncategorematic \emph{actual}
infinity for numbers and infinitesimals, following Ishiguro,
represents a conflation of the physical and the mathematical realms,
is unsupported by evidence, and must be rejected, as already indicated
by Antognazza (whose opinion Arthur quotes but fails to follow).

\subsection{Arthur's primary sources}
\label{s44}

In an attempt to provide a primary source for his application of
`syncategorematic actuality' to numbers, Arthur presents Leibniz's
\emph{New Essays} to the effect that:
\begin{quote}
``It is perfectly correct to say that there is an infinity of things,
  i.e. that there are always more of them than can be specified.  But
  it is easy to demonstrate that there is no infinite number, nor any
  infinite line or other infinite quantity, if these are taken to be
  genuine wholes. The Scholastics were taking this view, or should
  have been doing so, when they allowed a syncategorematic infinite,
  as they called it, and not a categorematic one.'' (A VI 6, 157/NE
  157).%
\footnote{Leibniz as reported in Arthur \cite[p.\,161]{Ar18}.}
\end{quote}
Thus far, Leibniz.  Arthur then proceeds to claim the following:
\begin{quote}
This is exactly the same notion of the actual infinite in multitude as
occurs in his mathematics: the number of parts (or bodies) is actually
infinite in the sense that for any finite number, there are actually
(not merely potentially) more bodies than this.  It is also precisely
the wording Leibniz had used in the \emph{New Essays} to describe what
the Scholastics meant, or `should have' meant `when they allowed a
syncategorematic infinite': `an infinity of things, i.e. that there
are always more of them than can be specified'.
\footnote{Arthur \cite[p.\,166]{Ar18}.}
\end{quote}
Arthur concludes:
\begin{quote}
I take the \emph{implicit criticism} in the \emph{`should have meant'}
to be that many Scholastics had followed Aristotle in identifying the
syncategorematic infinite with the potential infinite, with the result
that they denied the actual infinite.%
\footnote{Ibid.; emphasis added.}
\end{quote}
What was the `implicit criticism' that Arthur attributes to Leibniz?
We gather that it is not as though \emph{New Essays} ever mentioned
`syncategorematic actuality' for numbers; nonetheless, Arthur
speculates that this is what Leibniz must have had in mind when he
wrote that the Scholastics `{should have meant}'.  Such evidence is
clearly inconclusive.

There is a brief ``supplementary study'' \cite{Le06} where Leibniz
seems to illustrate syncategorematic actuality by numbers, but it
turns out that the study was ``Written on a separate slip of paper and
crossed through by Leibniz'' \cite[p.\;409]{Le07}; furthermore, the
supplementary study was never sent to des Bosses.  It seems
ill-advised to use such a flimsy study as a basis so as to re-evaluate
(as Arthur attempts to do in \cite[Section~10.5]{Ar18}) Leibniz's
well-established position concerning actuality as applicable only to
actual things.

\subsection{Analysis of letter to Varignon}
\label{s43}

As noted in Section~\ref{f29}, Ishiguro seeks to present Leibniz's
mention of the syncategorematic infinite as evidence for her
interpretation.  Let us consider the Leibnizian passage in the letter
to Varignon mentioning the syncategorematic infinite.  Leibniz wrote:
\begin{quote}
Yet we must not imagine that this explanation debases the science of
the infinite and reduces it to fictions, for there always remains a
`syncategorematic' infinite, as the Scholastics say {\ldots}%
%
\footnote{Leibniz \cite[p.\;93]{Le02} as translated by Loemker
(p. 543).}
%
\end{quote}
Does Leibniz mention the syncategorematic infinite in reference to an
infinite \emph{multitude} or an infinite \emph{magnitude} (see
Section~\ref{s34})?  It suffices to consider the second clause of
Leibniz's sentence:
\begin{quote}
And it remains true {\ldots}~that $2=\frac11 + \frac12 + \frac14 +
\frac18 +\frac1{16} +\frac1{32}+\,\cdots$, which is an infinite series
containing all the fractions whose numerators are $1$ and whose
denominators are a geometric progression of powers of $2$, although
only ordinary numbers are used, and no infinitely small fraction, or
one whose denominator is an infinite number, ever occurs in it.%
%
\footnote{Ibid., as translated by Loemker (pp.\;543--544).}
\end{quote}
Leibniz points out that even though the \emph{multitude} of the terms
of the series is infinite, it remains true that no infinite fractions
or infinite denominators are involved in summing the series.  Thus,
the term \emph{syncategorematic} in this passage does not refer to an
infinite \emph{magnitude}.

In modern mathematics, an infinite \emph{multitude} is formalized
e.g., by a Cantorian \emph{cardinal},%
\footnote{Leibniz has a more general notion of an `infinite whole';
see note~\ref{f32}.}
whereas an infinite \emph{magnitude} is formalized by a nonstandard
integer (see Section~\ref{s33}).  In his letter to Varignon, Leibniz
proceeds to discuss the multitude (not magnitude) of terms in a
geometric series.  Leibniz would of course reject the idea that such a
multitude exists as an `infinite whole', and would view it, instead,
syncategorematically (one can take as many of the terms as one
wishes).

Here Leibniz reiterates the point he made in his correspondence with
Bernoulli four years earlier, namely that convergent infinite series
cannot be used as evidence for the existence of infinitesimals.  The
passage does not support the inference sought by Ishiguro.%
\footnote{For further analysis of Ishiguro's hypothesis, see
\cite[Section 7.1]{16a}.}

According to Arthur,
\begin{quote}
Ishiguro claims that Leibniz `maintained that one can have a rigorous
language of infinity and infinitesimal while at the same time
considering these expressions as syncategorematic (in the sense of the
Scholastics), i.e., regarding the words as not designating entities
but as being well defined in the proposition in which they occur'
(1990, 82), citing Leibniz's words \emph{to this effect} in his letter
to Varignon of 1702.%
\footnote{Arthur \cite[p.\,168, note 19]{Ar18}; emphasis added.}
\end{quote}
Arthur claims that Leibniz's letter to Varignon corroborates
Ishiguro's hypothesis.  However, as we showed in this section,
Leibniz's words in his letter to Varignon were not at all to
Ishiguro's ``effect''.

\section{Arthur and Rabouin's erroneous reading of Robinson}
\label{s5}

Thus far, we have not mentioned either Abraham Robinson or his theory
NSA, except via direct quotations from Arthur and Rabouin.  The reason
is that AR's argument is more sweeping than merely rejecting NSA as a
historiographic tool; they reject any attribution of genuine
infinitesimals to the Leibnizian calculus (and consequently criticize
the work of Bos; see further in Section~\ref{s9}), making the question
of the applicability of NSA a secondary issue.

However, it would be remiss not to address AR's misinterpretation of
Robinson's discussion concerning the Leibnizian principle ``rules of
the finite succeed in the infinite, etc.'' as mentioned in their
note~4.%
\footnote{Arthur and Rabouin \cite[p.\;27, note 4]{Ar24}.}
We analyze the matter in Section~\ref{s51}.

\subsection{When do the rules of the finite succeed?}
\label{s51}

In their comments on Robinson's historical chapter, AR focus on the
issue concerning the precise nature of the ``principle'' that Leibniz
is referring to, and claim that Robinson misidentified the principle.
Regardless of the merits of their case, this issue is of secondary
importance since the main thrust of Robinson's comment lies elsewhere:
namely, the analogy with the transfer principle.  The Leibnizian
principle ``rules of the finite succeed in the infinite'' was not at
all obvious to mathematicians at the time, and in fact was the main
point of Leibniz's disagreement with Nieuwentijt.  Namely, Nieuwentijt
claimed that the square of an infinitesimal is precisely zero; see
e.g., Vermij \cite{Ve89}.  Nieuwentijt's stance evidently contradicts
Leibniz's comment concerning the rules of the finite (and similarly,
it would contradict the transfer principle).

AR's footnote 4 starts with a quotation from Robinson:
\begin{quote}
``The same principle [i.e.~Principle of Sufficient Reason]%
\footnote{The bracketed comment was added to Robinson's passage by
Arthur and Rabouin.}
is supposed to serve as justification for the assumption that the
infinitely small and infinitely large numbers obey the same laws as
ordinary (real) numbers.  Elsewhere, Leibniz appeals to his principle
of continuity to justify this assumption.''
\end{quote}
This passage appears on page 263 in Robinson \cite{Ro66}.  AR go on to
comment on Robinson's passage on page 263 as follows:
\begin{quote}
[Robinson's] first sentence rests on a misunderstanding of the letter
to Varignon from February 1702 in which Leibniz claims that the
``rules of the finite succeeds [sic] in the infinite as if there were
atoms'' and that, ``vice versa, the rules of the infinite succeeds
[sic] in the finite as if there were metaphysical infinitely smalls''
(although it is not the case that there are such entities in reality).
He then concludes: ``c'est parce que tout se gouverne par raison''.
But here, as in many other similar passages, the recourse to the
principle of reason is not a way to justify the transfer between
finite and infinite, but to explain why there are no atoms or true
infinitesimals in reality.%
\footnote{Arthur and Rabouin \cite[note~4]{Ar24}.}
\end{quote}
Here AR miss the point.  Robinson argued that his transfer principle
can be thought of as a formalisation of heuristic principles found in
Leibniz, such as his idea that ``the rules of the finite succeed in
the infinite {\ldots}~and vice versa, the rules of the infinite
succeed in the finite, etc.'' as we elaborate in Section~\ref{s52}.

\subsection{Transfer principle}
\label{s52}

To give an elementary illustration, the transfer principle entails
that a formula true for all real inputs would necessarily remain true
also for all hyperreal inputs.  For instance, consider the formula
expressing the idea that if a real number is nonzero, then its square
will also be nonzero, or in symbols
\[
\forall x\in \R\, (x\not=0 \to x^2\not=0).
\]
By transfer, the square of any nonzero infinitesimal must similarly be
nonzero.  To give another example, consider the trigonometric identity
$\sin^2x+\cos^2x=1$ for all $x\in\R$.  The trigonometric functions
possessing natural extensions to the hyperreal domain, the transfer
principle would therefore assert that the formula remains true for all
infinitesimal or infinite~$x$ (as well as all other hyperreal~$x$).

In their objection, AR miss the point of Robinson's insight.  Which
philosophical principle precisely Leibniz referred to in his letter to
Varignon is a secondary issue, the main issue being the parallelism
between Leibniz's idea of finite rules succeeding in the infinite, on
the one hand, and Robinson's transfer principle, on the other.  AR's
claim that Robinson misunderstood the letter to Varignon is based on
their own misunderstanding of the main thrust Robinson's comment on
page~263.  Robinson is even clearer on page~266:
\begin{quote}
Leibniz did say, in one of the passages quoted above, that what
succeeds for the finite numbers succeeds also for the infinite numbers
and vice versa, and this is remarkably close to our \emph{transfer} of
statements from~$\R$ to~$\astr$ and in the opposite direction.%
\footnote{Robinson \cite[p.\;266]{Ro66}.}
\end{quote}

Leibniz's heuristic principle (the Law of Continuity) was not a
self-evidence at the time, and in fact lay at the heart of his
argument against Nieuwentijt's approach (see \cite{Le95b}).
Nieuwentijt postulated a single infinite number, formed an
infinitesimal by passing to its reciprocal, and postulated that the
square of the infinitesimal is zero.  Nieuwentijt's infinitesimal
obviously does not adhere to the rule ``whatever succeeds for the
finite, succeeds also for the infinite'' (namely, it contradicts the
rule that a nonzero number must have a nonzero square).  Parallels
between the procedures of the Leibnizian calculus and those of modern
infinitesimal analysis are explored in \cite{24b}.

\subsection{Ehrlich on non-Archimedean theories}

In their discussion of ``non-Archimedean theories of the continuum in
general,'' AR claim that ``there were already \emph{good exemplars}
before the intervention of Robinson, as Philip Ehrlich has argued in
detail.''%
\footnote{Arthur and Rabouin \cite[p.\;28]{Ar24}; emphasis added.}

AR's claim amounts to both a slight toward Robinson's theory of
infinitesimals which was the first theory truly applicable in calculus
and analysis broadly understood (unlike the earlier theories), and a
misrepresentation of the position of Ehrlich, who wrote clearly:
\begin{quote}
After mathematicians had been taught for decades that a consistent
theory of the calculus based on infinitesimals was impossible, Abraham
Robinson was certainly swimming against the tide when he proved
otherwise.  \cite{Eh22}
\end{quote}
Ehrlich did present a detailed analysis of non-Archimedean theories
preceding Robinson's in \cite{Eh06}.  However, such theories failed
the Klein--Fraenkel criterion of applicability in broader mathematics,
as argued in~\cite{18i}.  More specifically, such theories did not
possess a sufficiently strong version of the transfer principle
discussed in Section~\ref{s52}.

\subsection
{Leibniz, Nelson, and inassignables}
\label{s54b}

Building upon Robinson's work in NSA in the 1960s \cite{Ro66}, Karel
Hrbacek \cite{Hr78} and Edward Nelson \cite{Ne77} developed axiomatic
approaches to NSA in the 1970s.  AR claim that Nelson's approach bears
resemblance to an early period (1669) when Leibniz endorsed actual
infinitesimals:
\begin{quote}
``[O]ne of us [namely, Arthur] has argued that prior to 1676 Leibniz
  did treat indivisibles as actual but unassignable elements of the
  continuum in three more or less distinguishable phases of his
  thought, in one of which they bear some resemblance to the
  non-Archimedean infinitesimals of Edward Nelson’s Internal Set
  Theory; see Arthur~(2009).''%
\footnote{Arthur and Rabouin \cite[p.\;28, note~13]{Ar24}.}
\end{quote}
The reference is to Arthur's 2009 article,%
\footnote{Arthur \cite[pp.\,15--16]{Ar09}.}
where an analogy with Nelson's work is applied to
\begin{quotation}
  ``Phase 1 (1669): Unassignable gaps in the continuum.''%
  \footnote{Op.~cit., p.\,39.}
\end{quotation}
This phase in Leibniz's development was analyzed in detail
by Ugaglia and Katz.%
\footnote{Ugaglia and Katz \cite[Section~3.2.1]{24e}.}
During this period,
\begin{quote}
[E]xtension is continuous, and continuity means for Leibniz
potentially infinite divisibility, that is syncategorematicity, as it
was for Aristotle.  The terminology itself of unassignability and
infinitesimals is absent in this phase:%
\footnote{Two years later, Leibniz will use \emph{inassignabilis} in
his \emph{Theoria motus abstracti} \cite[p.\;265]{Le71}.  This text
belongs to Phase II of Leibniz's development.}
`infinitely small' just means `potentially infinitely divisible.'%
\footnote{Ugaglia and Katz \cite[p.\,13]{24e}.}
\end{quote}
We will analyze this period in more detail in Section~\ref{s55}.

\subsection{Drift away from syncategoremata}
\label{s55}

In this early period, Leibniz endorsed the traditional syncategoremata
that he learned from the scholastics, which resulted in a
mathematically inconsistent approach.  There are no unassignable (or
inassignable) infinitesimals during this period.  Thus Arthur's
comparison of Leibniz's writings from 1669 and Nelson's theory is
unfounded.%
\footnote{Arthur seems to be aware of this difficulty when he writes
``The textual basis for Leibniz’s first theory of the continuum is
slender, and has to be pieced together from unpublished fragments and
his own later testimony'' \cite[p.\,14]{Ar09}.  We argue that Arthur
pieced it together incorrectly.}
Beeley writes:
\begin{quote}
Having little background in mathematics at that time [i.e., in
  1669--71], Leibniz had evidently gathered all he knew about
Cavalieri and Wallis and the criticisms which had been directed
against their respective methods from the polemical writings of Thomas
Hobbes.%
\footnote{Beeley \cite[pp.\;44--45]{Be08}.}
\end{quote}
As Leibniz's knowledge of mathematics grew, Leibniz abandoned the
earlier approaches (including syncategoremata) in favor of his
fictionalist approach to genuine infinitesimals.

It is instructive to ponder the implications of AR's discussion of
Nelson in connection with Leibniz's alleged period of ``Unassignable
gaps in the continuum.''  Such discussions, found in both Arthur
(2009) and Arthur and Rabouin (2024), appear to acknowledge that the
presence of unassignables (or inassignables) would signal
non-Archimedean behavior that ``bears some resemblance'' to Nelson's
theory, as already mentioned in Section~\ref{s4}.  AR's insights here
could be usefully extended to understanding Leibniz's mature views --
a route regrettably not taken by~AR.\, Curiously, Arthur does envision
comparing Leibniz's mature views with a modern theory called Smooth
Infinitesimal Analysis (SIA) in \cite{Ar13} (with the unsurprising
conclusion that they are not similar because SIA has \emph{nilpotent}
infinitesimals); see~\cite{21a} for an analysis.

\section{Hornangles \emph{vs} inassignables}

Hornangles (also known as angles of contingence, or angles of contact)
were discussed in Euclid and debated in the 16--17th centuries; see
Maier\`u \cite{Ma90}.  AR's discussion of hornangles overlooks the
fact that they are a manifestation of non-Archimedean behavior, as
acknowledged e.g., by Sonar, who writes:
\begin{quote}
[A]lready Eudoxus knew that also other number systems -- so-called
non-Archimedean number systems -- were conceivable {\ldots} Such a
system of quantities which was already known to the Greeks were
cornicular angles or horn angles.%
\footnote{Sonar \cite[p.\;40]{So21}.}
\end{quote}
Surprisingly, AR open their section 7 on hornangles with yet another
conflation of non-Archimedean phenomena and `nonstandard quantities',
as follows: ``it might be useful to analyse a last example regularly
mentioned in support of the existence of \emph{nonstandard quantities}
in Leibniz: horn angles.''%
\footnote{Arthur and Rabouin \cite[p.\;39]{Ar24}; emphasis added.}
AR fail to provide any quotation of an alleged attribution of a
``nonstandard'' attribute to hornangles, and with good reason: we are
not aware of any such attribution.%
\footnote{The reference \cite[Section~4.5]{24a} cited by AR in
footnote 99 on page~39 of \cite{Ar24} does contain a discussion of
hornangles, but certainly makes no such attribution, contrary to AR's
claim.}

\subsection{\emph{Euclidis Pr\^ota}: page 191}

AR's main source on hornangles is Leibniz's 1712 work \emph{Euclidis
Pr\^ota} \cite{Le12}.  They present the following passage on page 191
in GM:
\begin{quote}
``Indeed, it is not possible to assign a line which would be to
  another line as is the angle of contingency~$DCB$ at the angle of
  contingency~$DCF$, which we already showed to be possible for
  rectilinear angles.''%
\footnote{Leibniz as translated by Arthur and Rabouin \cite[p.\,39,
  note 101]{Ar24}.}
\end{quote}
AR go on to translate an additional passage from page 191:
\begin{quote}
``And therefore, they are from a completely different genre, and in
  relation to rectilinear angles, it [scil. the angle of contact]%
\footnote{The bracketed comment is in AR's text.}
  cannot even be considered as an infinitely small, which is always
  located between zero and an assignable.''%
\footnote{\label{f12}Leibniz as translated by Arthur and Rabouin
\cite[p.\;40]{Ar24}.  The original reads as follows: ``Ac proinde
plane est diversi generis, et respectu anguli rectilinei ne quidem ut
infinite parvus considerari potest, qui utique inter nullum et
assignabilem collocatur'' \cite{Ge50}, vol.\;V, p.\,191.}
\end{quote}
AR proceed to contest the conclusions of the analysis of Leibniz's use
of hornangles as developed by Katz et al.~\cite[Section~4.5]{24a}.  AR
rely mainly on these two quotations from page 191 in GM to develop
their rebuttal.  While they mention the year, 1712, of the Leibnizian
text \emph{Euclidis Pr\^ota}, they omit to mention the year of the
Leibnizian text that provided the basis for the analysis in
\cite{24a}.  That year is 1686 (Leibniz \cite{Le86}), namely no fewer
than 26 years prior to \emph{Euclidis Pr\^ota}.

\subsection{Leibniz on Peletier and Clavius}

AR similarly fail to quote a remarkable passage from page~191 in
\emph{Euclidis Pr\^ota} that clarifies Leibniz's position on both
hornangles and inassignables.  It turns out that Leibniz sides with
Peletier in his debate with Clavius (who endorsed hornangles),
declares that Euclid didn't really mean it when he described
hornangles as quantities, and says that Clavius' endorsement of
hornangles was a great embarrassment for geometers, and provided
fodder for Hobbes's ridicule of geometers:
\begin{quote}
``And so in this matter I agree with Peletarius against Clavius; and
  Euclid, when he said that the angle of contact is smaller than any
  straight line, spoke a little too loosely, understanding ``smaller''
  to mean something whose beginnings fall within the space of the
  former.  But it should not therefore be considered that he has
  attributed a perfect quantity to the angle of contact with respect
  to the straight line. {\ldots}
%
%
and Clavius erred greatly, when he denied this axiom, by which it is
asserted that what passes from one extreme to another, must also
equally pass through everything in between, and on this account gave
Thomas Hobbes an opportunity to insult geometers.''%
\footnote{Leibniz \cite[p.\,191]{Le12}; translation ours.}
\end{quote}

One is struck by the contrast between enthusiastic bashing of
hornangles in 1712, on the one hand, and their enthusiastic
endorsement 26 years earlier in his study of the degrees of contact
among curves in~\cite{Le86}, on the other.  The contrast indicates
that in the meantime Leibniz had had a change of heart concerning
hornangles.

\subsection{An infinitely small between zero and an assignable}
\label{s63}

But apart from the issue of hornangles, the most pertinent fact about
page~191 of \emph{Euclidis Pr\^ota} is that it contains a sentence
that directly refutes AR's insistence that Leibnizian infinitesimals
do not testify to a non-Archimedean phenomenon.  In fact, it is a
passage that was quoted by AR, who appear to have overlooked its
implications.  Indeed, when Leibniz writes that
\begin{quote}
``an infinitely small
is always located between zero and an assignable,''%
\footnote{See main text at note \ref{f12}.}
\end{quote}
he reveals that an infinitely small is inassignable and therefore
violates Euclid's Definition V.4, as Leibniz made clear elsewhere.%
\footnote{See note \ref{f5} for the references.}

A hornangle, while being smaller than every assignable rectilinear
angle, is not accorded the status of a genuine \emph{quantity} by
Leibniz in 1712.  Unlike a hornangle, an infinitesimal, similarly
smaller than every assignable quantity, is indeed accorded the status
of quantity as per page 191 of \emph{Euclidis Pr\^ota}.

There is a fundamental incoherence in AR's position with regard to
inassignables (or unassignables).  If inassignables are contradictory
and/or are to be understood syncategorematically (see
Section~\ref{f29}), why does Leibniz keep mentioning inassignables, as
for example when he describes the characteristic triangle as being of
inassignable size in his letter to Wallis, as analyzed by Beeley?%
\footnote{Leibniz \cite{Le98b}, quoted in Beeley \cite[p.\;48,
  note~36]{Be08}.  Beeley further mentions that ``[Leibniz], in a
letter to Wallis, {\ldots}~sets out reasons why the geometry of
indivisibles cannot strictly be reduced to the ancient method of
exhaustions, noting that the one operates with finite quantities, the
other with quantities incomparably smaller than the whole''
\cite[p.\;51]{Be08}.  Here Leibniz states clearly that his method does
\emph{not} operate merely with finite quantities, contrary to
Ishiguro's interpretation erdorsed by AR.}
Leibniz specifically objects to Wallis's habit of describing
infinitesimals as being exactly zero.  Leibniz points out that such an
approach introduces unnecessary obscurity into analysis.  Arguably,
dealing in contradictory entities would similarly introduce
unnecessary obscurity, as would be describing a
syncategorematically-inspired sequence of ordinary, assignable values
getting smaller and smaller as being \emph{inassignable}.%
\footnote{\label{f71}See main text in Section~\ref{s71} at
note~\ref{f74} for a related point concerning bounded infinities.}

\section{Proposition\;11 and \emph{infinita terminata}}
\label{s7}

In his \emph{De Quadratura Arithmetica}, Leibniz introduces the notion
of \emph{infinita terminata} (literally ``bounded infinites'') and
contrasts it with unbounded infinites, or infinite wholes; see
Section~\ref{s3}.

\subsection{Should Leibniz be taken literally?}
In this connection, AR write:
\begin{quote}
Leibniz’s strategy in prop.\;XI is to shortcut any reference to this
`unbounded' space (for which the notion of area would be problematic)
by introducing a parallel line to the asymptote which can be taken as
close as one wishes but will still intersect the given curve.%
%
~{\ldots}\;But that means that the line, which was introduced to
perform the quadrature, can be made as long as one wants (greater that
any given magnitude).%
\footnote{Arthur and Rabouin \cite[p.\;34]{Ar24}.}
\end{quote}
They go on to the literal meaning:
\begin{quote}
\emph{Taken literally (as a fixed and not as a variable object), it
hence describes a line which is at the same time bounded} (it
intersects the curve at a point and we can make it begin at the origin
of our frame of coordinates) \emph{and infinite} (it can be made
larger than any given quantity).%
\footnote{Op.~cit., pp.\;34--35; emphasis added.}
\end{quote}
Note that AR are not quoting Leibniz here but only paraphrasing him.
How do they mean exactly to say that ``Taken literally (as a fixed and
not as a variable object), it hence describes a line which is at the
same time bounded {\ldots}~and infinite''?  Did Leibniz mean it
literally, or did he mean to say that ``it can be made larger than any
given quantity" as AR claim?  These are apparently two different
things, in spite of AR's effort to blur the distinction.  Either
Leibniz meant it literally, or else what he had in mind was a
paraphrase along the lines ``can be made as long as one wants".

It is significant that AR don't quote Leibniz directly, because if one
examines what Leibniz actually wrote (see for instance \cite{23h}),
there is little doubt as to what he meant, as we discuss in
Section~\ref{s71}.

\subsection{Leibniz's infinitesimal $(\mu)\mu$}
\label{s71}
 
In Proposition 11, Leibniz proceeds as follows.

\begin{enumerate}
\item
He denotes a point at infinitesimal distance from the origin
by~$(\mu)$, where the parentheses are part of the notation
(whereas~$\mu$ itself denotes the origin).
\item
He considers a line parallel to the asymptote and at infinitesimal
distance from it (i.e., passing through the point~$(\mu)$ as above).
\item
He takes the corresponding point on the curve, whose ordinate is then
a specific fixed (and not a variable) infinite entity, namely an
\emph{infinitum terminatum}.
\end{enumerate}

AR's paraphrase of Leibniz rests upon the assumption that he must have
meant it syncategorematically (see Section~\ref{f29}).  In doing so,
they create an impression that a genuine bounded infinity would have
to be as contradictory as a quantity that is both fixed and variable.

AR's conclusion with regard to the Leibnizian \emph{bounded infinity}
is that ``the line, which was introduced to perform the quadrature,
can be made as long as one wants (greater that any given magnitude)''
as quoted above.  However, this begs the question concerning their
Ishiguro/exhaustion reading.  According to their explanation, there is
no difference between bounded infinity and an (arbitrarily large)
finite line segment.  If so, what would be the point for Leibniz to
introduce a paradoxical-sounding term (``bounded infinity''!) if he
simply meant `finite line segment'?  This amounts to a basic
incoherence of their intepretation.%
\footnote{\label{f74}See main text in Section~\ref{s63} at note
\ref{f71} for a related point concerning infinitesimals.}
%

%

\subsection{Three approaches to bounded infinity}

How does Leibniz deal with the contradictory nature of an infinite
whole such as an unbounded line?  He dissolves the contradiction by
saying that one can take a finite line bigger than any given, i.e., by
viewing infinity syncategorematically (\emph{not} as a whole).  How do
AR interpret bounded infinity?  Also by saying that one can take a
finite line bigger than any given.

It emerges that, on AR's reading, \emph{infinitum terminatum} is
indistinguishable from unbounded infinity.  Such an interpretation
however is incoherent because Leibniz explicitly contrasts them in the
Scholium to Proposition~11.%
\footnote{See note \ref{f23}.}

There are three interpretations of \emph{infinitum terminatum} in the
current literature:
\begin{enumerate}
\item
As argued in \cite{23h}, on Knobloch's reading it amounts to an ideal
perspective point at infinity.%
\footnote{Briefly, Knobloch analyzes theorem 45 in Leibniz's \emph{De
Quadratura Arithmetica}, and claims that the point at infinity
introduced there is an instance of the \emph{infinitum terminatum}.
However, that point at infinity is recognizably the ideal perspective
point at infinity known already to Kepler and Desargues, as argued in
\cite{23h}.}
\item
On AR's Ishiguro-inspired reading, it is indistinguishable from
unbounded infinity.
\item
We argue that it is an inverted infinitesimal and thus a genuine
mathematical entity not found in the physical realm.
\end{enumerate}
The weaknesses of interpretations (1) and (2) are analyzed
respectively in \cite{23h} and the present text.

\section{Question of use and question of existence}
\label{s8}

In interpreting historical infinitesimalists, it is important to
distinguish between mathematical practice, on the one hand, and
foundational ontology, on the other.  Modern analysis, whether of
Weierstrassian or Robinsonian variety, is typically built upon
set-theoretic foundations involving infinite sets as a \emph{sine qua
non} (but see Section~\ref{s33}).  Such foundations are far removed
from the world of the historical infinitesimalists such as Leibniz,
Euler, and Cauchy.  On the other hand, the \emph{procedures} of modern
analysis, whether of Weierstrassian or Robinsonian kind, often find
close analogs in the work of historical infinitesimalists.

It is obvious that a study of Leibnizian calculus must focus on the
procedures if it is to avoid the trap of presentism.  We can only
study the procedures Leibniz used in his calculus and geometry, and
certainly cannot attribute modern set-theoretic investigations to him.
Thus the issue of use should be clearly separated from the issue of
existence as a mathematical entity in a set-theoretic framework.

AR appear to set this distinction on its head when they claim that
\begin{quote}
[A]s long as one conflates questions of use and questions of
existence, it will always be possible to counterbalance any explicit
denial of the existence of infinitesimals by Leibniz by an appeal to
his practice as testifying to the contrary.  At a minimum, \emph{one
could argue} that there is tension between what he says and what he
does.%
\footnote{Arthur and Rabouin \cite[p.\;27]{Ar24}; emphasis added.}
\end{quote}

However, this is not at all what has been argued in recent Leibniz
scholarship on his genuine infinitesimals.  Leibniz clearly denied
that infinitesimals exist in the phenomenal realm, and they are
therefore fictional entities like negatives and imaginary roots.
Contrary to AR's claim, such a reading does not constitute a
``conflation of the question of use with that of existence of
infinitely small quantities.''  AR go on to claim that
\begin{quote}
Leibniz himself emphasized on many occasions that his practice should
remain neutral as regards the question of existence, and henceforth
that the use of infinitesimals should be interpretable in both ways
(i.e., either as existing entities or not).%
\footnote{Op.~cit., p.\;28.}
\end{quote}
Their parenthetical comment reveals a basic misunderstanding of the
Leibnizian position.  On many occasions, Leibniz emphasized that the
use of infinitesimals can be \emph{replaced} by equivalent (but more
laborious) procedures \`a la Archimedean exhaustion.  It is misleading
to characterize the existence of two distinct methods as treating
infinitesimals ``either as existing entities or not'' since by all
accounts, they do not exist in the phenomenal realm, according to
Leibniz.  The issue is whether or not to use them directly.  AR
conclude that
\begin{quote}
[One] needs to explain how it is that one can employ fictions of this
type (contradictory notions) in mathematical discourse without falling
into inconsistency.  Moreover, if one conflates questions of use and
questions of existence, this explanation will not alter the impression
that there is a \emph{tension} between Leibniz’s justification and his
practice.%
\footnote{Op.~cit., p.\;41; emphasis added.}
\end{quote}

If there is any `tension', it is created by AR's misconceived idea
that infinitesimals are `contradictory notions'.  The alleged
conflation between use and existence is a made-up problem, as
explained above.

\section{Two syncategorematists against Bos}
\label{s9}

The influential study of Leibnizian methodology by Bos \cite{Bo74} was
published in 1974.  Both Arthur and Rabouin are on record as
criticizing the study by Bos whenever his analysis creates tensions
with Ishiguro's 1990 reading (see Section~\ref{f29}).


In 2020, Rabouin and Arthur quote an example of Leibniz's application
of the Law of Continuity, and then proceed to criticize Bos in the
following terms:

\begin{quote}
This first basic example was not mentioned in Henk Bos’s seminal
study (Bos 1974–5), and this had the unfortunate consequence of hiding
the \emph{global structure} of the reasoning.%
\footnote{Rabouin and Arthur \cite[p.\;439]{Ra20}; emphasis added.}
\end{quote}
The alleged `global structure' is Ishiguro's syncategorematic reading.
In reference to an additional example, they lodge a further criticism:

\begin{quote}
Here again, the presentation by Bos, as thorough and interesting as
it is, had the unfortunate consequence of hiding some particularities
of the reasoning by neglecting some developments that were apparently
of no interest.%
\footnote{Ibid.}
\end{quote}
Related criticisms by Arthur appear in an earlier text dating from
2013.


\subsection{Arthur \emph{vs} Bos, 2013}

In 2013, Arthur claimed the following:

\begin{quote}
[I]t is \emph{not} the case that Leibniz has two methods, one
  committed to the existence of infinitesimals as `genuine
  mathematical entities', `fixed, but infinitely small', and the other
  Archimedean, treating the infinitely small as arbitrarily small
  finite lines, ones that could be made as small as desired.%
\footnote{Arthur \cite[p.\;561]{Ar13}; emphasis added.}
\end{quote}
Bos's position was precisely that Leibniz \emph{did} have two methods;
we concur.  Arthur continues criticizing Bos in footnote 5 on the same
page:

\begin{quote}
Bos, in his classic article on Leibniz’s differentials (1974–1975,
55), sees Leibniz as pursuing ``two different approaches to the
foundations of the calculus; one connected with the classical methods
of proof by ‘exhaustion’, the other in connection with a law of
continuity.''  The quotations given in my text are from Bos
(1974–1975), 12, 13.  See also Herbert Breger's similar
\emph{criticisms of Bos} on these points in his (2008, 195–197).%
\footnote{Op.~cit., note 5; emphasis added.}
\end{quote}
Arthur goes on to contrast his work with that of Bos:
\begin{quote}
My interpretation here contrasts with that of Bos, who reads Leibniz
as having two different interpretations of differentials, one in which
they are infinitely small elements of lines, and the other where they
stand for finite lines.%
\footnote{Op.~cit., p.\;566.}
\end{quote}
Here Bos summarized Leibniz's approach of using~$dx, dy$ for
infinitesimal differentials, and a separate notation~$(d)x, (d)y$ when
``ordinary'' values are substituted for~$dx,dy$.%
\footnote{In his 1974 article, Bos used~$\underline d x$ in place of
Leibniz's~$(d)x$.}
After outlining the reading by Bos, Arthur goes on to lodge an even
stronger claim of ``misinterpretation'':
\begin{quote}
But this, it seems to me, is to \emph{misinterpret} Leibniz’s
justification, {\ldots}%
\footnote{Arthur \cite[p.\;567]{Ar13}; emphasis added.}
\end{quote}

\subsection{Who is Arthur and Rabouin's quarrel with?}

It emerges from the analysis above that AR's quarrel is neither with
Robinson, nor Nelson, nor yet ``the mathematician Mikhail G. Katz with
various co-authors''%
\footnote{Arthur and Rabouin \cite[note 2]{Ar24}.}
but rather squarely with the historian Henk Bos, whose work on Leibniz
is described by Guicciardini as ``a classic in our field''
\cite[pp.\;40--41]{Gu24}.  Naturally, other historians have also
disagreed with AR's interpretation.  Thus, the 2021 article by
Esquisabel and Raffo Quintana explicitly rejects it in the following
terms:
\begin{enumerate}
\item
``[U]nlike the infinite number or the number of all numbers, for
  Leibniz infinitary concepts do not imply any contradiction, although
  they may imply paradoxical consequences.''  \cite[p.\;641]{Es21}
\item
``[W]e disagree with the reasons [Rabouin and Arthur] gave for the
  Leibnizian rejection of the existence of infinitesimals, and in our
  opinion the texts they refer to in order to support their
  interpretation are not convincing.  Since we argue that Leibniz did
  not consider the concept of infinitesimal as
  \emph{self-contradictory}, we try to provide an alternative
  conception of impossibility.''\;\cite[p.\;620]{Es21}
\end{enumerate}
Antognazza's objections to Arthur's interpretation were analyzed in
Section~\ref{s21b}.  Further objections were detailed in Eklund
\cite{Ek20}.

\subsection{Jesseph vs the syncategorematic interpretation}
\label{jess}

Jesseph sees Leibniz as having two ways of interpreting the calculus:
\begin{quote}
Rather than postulate an ambiguity in Leibniz’ writings, we can take
him at his word: he is convinced that infinitesimal magnitudes are
eliminable, and the reason for this conviction is the thoroughly
commonsensical belief that the truth of his mathematical results
should not depend too heavily upon the resolution of metaphysical
problems.  In the Leibnizian scheme, true mathematical principles will
be found acceptable on any resolution of the metaphysical problems of
the infinite. Thus, Leibniz' concern with matters of rigor leads him
to propound a very strong thesis indeed, namely no matter how the
symbols ``dx'' and ``dy'' are interpreted, the basic procedures of the
calculus can be vindicated.%
\footnote{Jesseph \cite[p.\;243]{Je89}.}
\end{quote}
Jesseph concludes:
\begin{quote}
Such vindication could take the form of a
new science of infinity, or it could be carried out along classical
lines, but in either case the new methods will be found completely
secure.%
\footnote{Ibid.}
\end{quote}
Accordingly, infinitesimals are not merely an abbreviated way of
speaking (as per to Ishiguro) and certainly not contradictory (as per
AR), but rather well-founded fictions.  In his 1998 text on page~35,
Jesseph explains the difference between well-founded and ill-founded
fictions, and concludes:
\begin{quote}
Infinitesimal magnitudes {\ldots}~are nevertheless well-founded
fictions precisely because the realm of mathematical objects is
structured just as it would be if the infinitesimal existed.%
\footnote{Jesseph \cite[pp.\;35--36]{Je98}.}
\end{quote}
In these passages, Jesseph clearly distances himself from the
syncategorematic interpretation of Leibniz's infinitesimals.  It is
therefore surprising that he should have consented to have his name
added to the author list of the broadside \cite{Ar22}.

To summarize, we list some of the Leibniz scholars who disagree with
the reading by Ishiguro, Arthur, and Rabouin: Antognazza, Bos, Eklund,
Esquisabel and Quintana, Jesseph.

\section{Arthur vs Rabouin: Limit processes in 2013}
\label{s10}

In 2013, a pair of texts appeared that offer an interesting
comparison.  One of them, ``Leibniz's syncategorematic
infinitesimals'' by Arthur, claims that (syncategorematic)
infinitesimals involve a `limiting process'.  Another, by Rabouin,
claims that, unlike the Newtonian tradition, there were no `limit
processes' in the Leibnizian calculus.  Thus, Arthur wrote:
\begin{quote}
The role of~$(d)y$,~$(d)v$ and~$(d)x$ is to be ``finite surrogates'':
they remain assignable quantities throughout the \emph{limiting
\mbox{process}} as~$dx$ becomes arbitrarily small, and since~$(d)y$
and $(d)x$ have an assignable ratio equal to~$dy:dx$ even in that
\emph{limit} (likewise for~$(d)v$ and~$(d)x$), differential equations
involving terms in~$dy$,~$dv$ and~$dx$ to the same order are
interpretable under the fiction where these quantities stand for
infinitely small quantities.%
\footnote{Arthur \cite[p.\;567]{Ar13}; emphasis added.}
\end{quote}
Arthur's is a clear statement to the effect that the syncategorematic
interpretation involves a `limiting process'.  Meanwhile, Rabouin
wrote in the same year:
\begin{quote}
A profound evolution of the concept involved in differential calculus,
already emphasized by D'Alembert, was the progressive introduction of
a new interpretation of its basic constituents in terms of ``limit
processes'' -- an interpretation that is \emph{absent from the
practice of Leibniz} and his first followers, but present in the
Newtonian tradition.%
\footnote{Rabouin \cite[p.\;24]{Ra13}; emphasis added.  By 2017, the
passage was modified by deleting the phrase ``absent from the practice
of Leibniz'' as follows: ``A profound evolution of the concept
involved in differential calculus, already emphasized by Jean le Rond
d'Alembert, was the progressive introduction of a new interpretation
of its basic constituents in terms of limit processes--an
interpretation that is put forward in the Newtonian tradition.  The
Leibnizian calculus developed quite rapidly as a mixture of features
derived from Leibniz's and Newton's techniques, etc.''
\cite[p.\;213]{Ra17}.  The 2013 version was the more accurate one.}
\end{quote}

We concur with Rabouin's view as expressed in 2013 that in Leibniz's
mature view of the calculus, there were no limit processes -- and
therefore no syncategorematics as per Arthur 2013.

If the interpretation in terms of ``limit processes'' is absent from
the Leibnizian calculus as per Rabouin 2013, then infinitesimals are
not to be interpreted syncategorematically (as per Arthur 2013), and
therefore they are entities exhibiting non-Archimedean behavior as per
Rabouin 2013.  Since there is no trace of any endorsement of an NSA
interpretation of Leibniz in Rabouin 2013, it emerges that Rabouin was
capable of envisioning non-Archimedean phenomena in Leibniz without
endorsing any concomitant Leibniz--Robinson linkage.  Somehow the
distinction between genuine infinitesimals and those of Robinsonian
analysis has been erased by the time of writing Rabouin's 2024 text
with Arthur.

\section{Conclusion}

We have analyzed the recent publication by Arthur and Rabouin in
\emph{Historia Mathematica} on Leibniz's infinitesimals, and
identified many errors stemming in a number of cases from flawed
analyses in their earlier publications.  Their thesis is in fact far
more radical than what their title and abstract may suggest.  Arthur
and Rabouin are not merely rejecting interpretations of the Leibnizian
calculus based on Abraham Robinson's theory.  Their claim goes beyond
the issue of the historiographic applicability of nonstandard
analysis, and amounts to a sweeping claim concerning the Leibnizian
calculus itself.  Their claim consists in rejecting any attribution of
genuine infinitesimals to the Leibnizian calculus.

It is no accident that their article fails to cite the seminal study
of Leibnizian methodology by Bos \cite{Bo74}, described as ``a classic
in our field'' by Guicciardini \cite{Gu24}.  Bos did acknowledge that
Leibniz used genuine infinitesimals (and there is no trace of
syncategoremata in his study).  In sum, Arthur and Rabouin's claim is
unviable.


\section{Acknowledgments} We are grateful to Ali Enayat,
Karel Hrbacek, David Sherry, Monica Ugaglia, and the anonymous referee
for helpful suggestions.  The influence of Hilton Kramer (1928--2012)
is obvious.

\end{document}